\documentclass[dvips,11pt,a4paper, oneside]{amsart}
\usepackage{amssymb}
\usepackage[latin2]{inputenc}
\usepackage{graphicx}
\usepackage{indentfirst}

\DeclareMathOperator{\diam}{diam}

\newtheorem{Def}{Definition}
\newtheorem{Lem}{Lemma}
\newtheorem{Prop}{Proposition}
\newtheorem{Thm}{Theorem}
\newtheorem{Cor}{Corollary}
\newtheorem{Rem}{Remark}
\newenvironment{Pf}{ Proof.}{\(\square\)}

\title[Generalized conic functions ...]{Generalized conic functions of hv-convex planar sets: continuity properties and relations to X-rays}
\author{Csaba Vincze and \'Abris Nagy}
\address{Csaba Vincze}
\address{Inst. of Math., Univ. of Debrecen, \\
H-4010 Debrecen, P.O.Box 12 \\
Hungary.}
\email{csvincze@science.unideb.hu}
\address{\'Abris Nagy}
\address{Inst. of Math., MTA-DE Research Group `Equations Functions and Curves',}
\address{Hungarian Academy of Sciences.}
\address{Inst. of Math., Univ. of Debrecen, H-4010 Debrecen, P.O.Box 12 Hungary.}
\email{abris.nagy@science.unideb.hu}

\begin{document}
\begin{abstract}
In the paper we investigate the continuity properties of the mapping $\Phi$ which sends any non-empty compact connected hv-convex planar set $K$ to the associated generalized conic function $f_K$. The function $f_K$ measures the average taxicab distance of the points in the plane from the focal set $K$ by integration. The main area of the applications is the geometric tomography because $f_K$ involves the coordinate X-rays' information as second order partial derivatives \cite{NV3}. We prove that the Hausdorff-convergence implies the convergence of the conic functions with respect to both the supremum-norm and the $L_1$-norm provided that we restrict the domain to 
the collection of non-empty compact connected hv-convex planar sets contained in a fixed box (reference set) with parallel sides to the coordinate axes. We also have that $\Phi^{-1}$ is upper semi-continuous as a set-valued mapping. The upper semi-continuity establishes an approximating process in the sense that if $f_L$ is close to $f_K$ then $L$ must be close to an element $K'$ such that $f_{K}=f_{K'}$. Therefore $K$ and $K'$ have the same coordinate X-rays almost everywhere. Lower semi-continuity is usually related to the existence of continuous selections. If a set-valued mapping is both upper and lower semi-continuous at a point of its domain it is called continuous. The last section of the paper is devoted to the case of non-empty compact convex planar sets. We show that the class of convex bodies that are determined by their coordinate X-rays coincides with the family of convex bodies $K$ for which $f_K$ is a point of lower semi-continuity for $\Phi^{-1}$.
\end{abstract}
\maketitle
\section{Introduction}

The idea motivating our investigations is the application of generalized conics' theory \cite{NV1}, \cite{NV2} and \cite{NV3} in geometric tomography. Let $K$ be a compact planar set in the Euclidean plane $\mathbb{R}^2$, $p\geq 1$ and consider the distance function 
$$d_p\big((x,y),(\alpha,\beta)\big)=\sqrt[p\ ]{|x-\alpha|^p+|y-\beta|^p}$$
induced by the p - norm. Pairs of the form $(x,y)$, $(\alpha,\beta)$ and $(c_1,c_2)$ denote elements of $\mathbb{R}^2$. Let us define the sets
	\[x<_1 K:=\{(\alpha,\beta)\in K\ | \ x<\alpha\},\ \ K<_1 x:=\{(\alpha,\beta)\in K\ |\  \alpha<x\},
\]
	\[y<_2 K:=\{(\alpha,\beta)\in K\ | \ y<\beta\},\ \ K<_2 y:=\{(\alpha,\beta)\in K\ | \ \beta<y\},
\]
	\[x=_1 K:=\{(\alpha,\beta)\in K\ |\  \alpha=x\},\ \ y=_2 K:=\{(\alpha,\beta)\in K\ | \ \beta=y\},
\]where the index refers to the usual ordering of the coordinates. A planar set $K$ is said to be hv-convex if the sections $x=_1 K$ and $y=_2 K$ are convex sets for all $x,y\in\mathbb{R}$.

\begin{Def}The \emph{X-ray functions into the coordinate directions} are
	\[Y_K(x):=\lambda_1(x=_1 K)\ \ \textrm{and} \ \ X_K(y):=\lambda_1(y=_2 K),
\]where $\lambda_1$ denotes the one-dimensional Lebesgue measure.
\end{Def}

\noindent
X-ray functions (especially coordinate X-rays) are typical objects in geometric tomography \cite{Gardner}. We start from a compact set $K$ in the plane to construct a convex function carrying the information of coordinate X-rays as second order derivatives.
\begin{Def} \emph{\cite{NV3}}
\emph{The generalized conic function} $f_K$ associated to $K$ is defined by the formula
	\[f_K(x,y):=\intop\limits_K d_1\big((x,y),(\alpha,\beta)\big)\, d\alpha d\beta.
\]
The levels of this function are called \emph{generalized conics} with $K$ as the \emph{focal set}.
\end{Def}

Using that the 1 - norm is decomposable the generalized conic function $f_K$ can be expressed in terms of the coordinate X-rays as follows:
	\begin{equation}
	\label{key1}
	f_K(x,y)=\intop\limits_{-\infty}^{\infty} |x-\alpha| Y_K(\alpha)\ d\alpha+\intop\limits_{-\infty}^{\infty} |y-\beta| X_K(\beta)\ d\beta;
\end{equation}
moreover
	\[\frac{\partial f_K}{\partial x}(x,y)=\lambda_2(K<_1 x)- \lambda_2(x<_1 K)
\]and, in a similar way,
	\[\frac{\partial f_K}{\partial y}(x,y)=\lambda_2(K<_2 y)-\lambda_2(y<_2 K),
\]where $\lambda_2$ denotes the Lebesgue measure on $\mathbb{R}^2$. By the Cavalieri's principle
	\[\lambda_2(K<_1 x)=\intop \limits_{-\infty}^x Y_K(s)\, ds\ \ \ \textrm{and}\ \ \ \ \lambda_2(x<_1 K)=\intop \limits_x^{\infty} Y_K(s)\, ds,
\]
	\[\lambda_2(K<_2 y)=\intop \limits_{-\infty}^y X_K(t)\, dt \ \ \ \textrm{and}\ \ \ \ \lambda_2(y<_2K)=\intop \limits_y^{\infty} X_K(t)\, dt.
\]Lebesgue differentiation theorem leads us to
	\begin{equation}
	\label{key2}
	\frac{\partial^2 f_K}{\partial x \partial x}(x,y)=2 Y_K(x)\ \ \textrm{and}\ \ \ \frac{\partial^2 f_K}{\partial y \partial y}(x,y)=2 X_K(y)
\end{equation}
except on a set of measure zero. Equations (\ref{key1}) and (\ref{key2}) show that $f_K=f_L$ if and only if $K$ and $L$ have the same coordinate X-rays almost everywhere. 

\begin{Rem} If $\lambda_2(K)>0$ then the weighted function
$$F_K:=\frac{1}{\lambda_2(K)}f_K$$
can be also introduced. We have that $F_K=F_L$ if and only if the coordinate X-rays are proportional to each other \cite{NV3}. 
\end{Rem}

Our main result is that 
the Hausdorff-convergence implies the convergence of the conic functions with respect to both the supremum-norm and the $L_1$-norm provided that we restrict the domain to 
the collection of non-empty compact connected hv-convex planar sets contained in a fixed box (reference set) with parallel sides to the coordinate axes. We also have that $\Phi^{-1}$ is upper semi-continuous as a set-valued mapping. The upper semi-continuity establishes an approximating process in the sense that if $f_L$ is close to $f_K$ then $L$ must be close to an element $K'$ such that $f_{K}=f_{K'}$. Therefore $K$ and $K'$ have the same coordinate X-rays almost everywhere. R. Gardner and M. Kiderlen \cite{KG} presented an algorithm for reconstructing convex bodies from noisy X-ray measurements with a full proof of convergence in 2007. In this sense the reconstruction means to give the unknown set as a limit of a convergent sequence. The algorithm uses four directions which is related to the minimal number of directions for all convex bodies to be determined by their X-rays in these directions. Our approach means an alternative way for the reconstruction/approximation. The method uses only two directions (coordinate X-rays) and we can apply the theory to the wider class of non-empty compact connected hv-convex planar sets: let $f_K$ be the input data and consider the optimization problem
$$\textrm{minimize} \ \|f_L-f_K\|\ \ \textrm{subject to}\ \  L\in \mathcal{H},$$ 
where the set $\mathcal{H}$ is the subcollection of non-empty compact connected hv-convex sets which are constituted by the subrectangles belonging to a partition of the reference set \cite{NV4}. It is typically a rectangle $B$ with parallel sides to the coordinate axes such that $K\subset B$.

\section{General observations}

In what follows $\mathcal{F}(\mathbb{R}^2)$ denotes the metric space of non-empty bounded and closed (i.e. compact) subsets in the plane equipped with the Hausdorff metric \cite{Coppel}. The outer parallel body $K^{\varepsilon}$ is the union of all closed Euclidean balls centered at the points of $K\in \mathcal{F}(\mathbb{R}^2)$ with radius $\varepsilon>0$. The Hausdorff distance between $K$ and $L\in \mathcal{F}(\mathbb{R}^2)$ is given by the formula
	\[H(K,L):=\inf \{\varepsilon >0\ | \ K\subset L^{\varepsilon}\ \textrm{\emph{and}}\ L\subset K^{\varepsilon}\}.
\]
Let $\textrm{pr}_1$ and $\textrm{pr}_2$ be the orthogonal projections onto the coordinate axes and consider a rectangle $B$ with parallel sides to the coordinate axes. The level set of $\textrm{pr}_1\times \textrm{pr}_2$ belonging to $B$ is defined as
$$\mathcal{L}_B=\{ L \in \mathcal{F}(\mathbb{R}^2)\ | \ \textrm{pr}_1(L)\times \textrm{pr}_2(L)=B\}.$$
We also introduce a kind of sublevel set 
$$\mathcal{M}_B=\bigcup_{B^*\subset B} \mathcal{L}_{B^*}$$
with respect to the partial ordering induced by the inclusion, where the union is taken with respect to all rectangles with parallel sides to the coordinate axes and contained in $B$.

\begin{Prop}\label{prop2}
Both $\mathcal{L}_B$ and $\mathcal{M}_B$ are compact with respect to the Hausdorff metric.
\end{Prop}

\begin{Pf}
Since $\mathbb{R}^2$ is complete, the space $\mathcal{F}(\mathbb{R}^2)$ equipped with the Hausdorff metric is also complete, see [1]. By Hausdorff's theorem any closed and totally bounded subset in a complete metric space is compact. That $\ \mathcal{L}_B$ (or $\mathcal{M}_B$) is totally bounded follows from the well-known version of Blaschke's selection theorem for compact sets \cite{Schneider} [Theorem 1.8.4]. The closedness can be concluded from the continuity of the mapping
$$L \mapsto \textrm{pr}_1(L)\times \textrm{pr}_2(L),$$
where $\textrm{pr}_1$ and $\textrm{pr}_2$ denote the orthogonal projections onto the coordinate axes. The orthogonal projections are obviously continuous mappings with respect to the Hausdorff metric because the projection of any closed ball is a lower dimensional closed ball. In other words the projected parallel body is just the parallel body (with the same radius) of the projected set. 
\end{Pf}

\begin{Prop}
Both $\mathcal{L}_B$ and $\mathcal{M}_B$ are convex in the sense that $L_1\in \mathcal{L}_B$/$\mathcal{M}_B$ and $L_2\in \mathcal{L}_B$/$\mathcal{M}_B$ implies that $tL_1+(1-t)L_2\in \mathcal{L}_B$/$\mathcal{M}_B$
for all $\ 0\leq t\leq 1$.
\end{Prop}

\begin{Pf} Since the projections preserve the convex combination of the elements the statement follows directly from the definition of $\mathcal{L}_B$/$\mathcal{M}_B$.
\end{Pf}

\begin{Thm} The mapping $\Phi\colon L\in \ \mathcal{L}_B\to f_L$ is concave in the sense that
for any $(x,y)\in \mathbb{R}^2$
	\[f_{tL_1+(1-t)L_2}(x,y)\geq tf_{L_1}(x,y)+(1-t)f_{L_2}(x,y),
\]where $L_1$, $L_2\in \mathcal{L}_B$ and $\ 0\leq t \leq 1$.
\end{Thm}

\begin{Pf}
Let $(\alpha, \beta)\in B$ be a fixed point. Then
	\[\big(\alpha=_1 (tL_1+(1-t)L_2) \big)\supset \big(t\alpha=_1 tL_1 \big)+\big((1-t)\alpha=_1 (1-t)L_2 \big),
\]
	\[\big(\beta=_2 (tL_1+(1-t)L_2) \big)\supset \big(t\beta=_2 tL_1 \big)+\big((1-t)\beta=_2 (1-t)L_2 \big).
\]Since none of the sets $\alpha=_1 L_i$ ($i=1,2$) are empty the numbers
	\[c_1:=\inf \{y\ | \ (t\alpha,y)\in tL_1\}=t\inf \{y\ |\ (\alpha,y)\in L_1\},
\]
	\[d_1:=\sup \{y\ | \ (t\alpha,y)\in tL_1\}=t\sup \{y\ |\ (\alpha,y)\in L_1\}
\]and
	\[c_2:=\inf \{y\ | \ ((1-t)\alpha,y)\in (1-t)L_2\}=(1-t)\inf \{y\ |\ (\alpha,y)\in L_2\},
\]
	\[d_2:=\sup \{y\ | \ ((1-t)\alpha,y)\in (1-t)L_2\}=(1-t)\sup \{y\ |\ (\alpha,y)\in L_2\}
\]are well-defined. Consider the sets
	\[P:=(t\alpha,c_1)+(1-t)(\alpha=_1 L_2)\subset \textrm{conv}\ \big\{(\alpha,c_1+c_2),(\alpha,c_1+d_2)\big\}
\]and
	\[Q:=((1-t)\alpha,d_2)+t(\alpha=_1 L_1)\subset \textrm{conv}\ \big\{(\alpha,d_2+c_1),(\alpha,d_2+d_1)\big\}.
\]Since they are contained in adjacent segments and
	\[P\cup Q\subset \big(\alpha =_1 tL_1+(1-t)L_2 \big)
\]it follows that
	\[Y_{t L_1+(1-t) L_2}(\alpha)\geq \lambda(P\cup Q)=\lambda_1 (P) + \lambda_1 (Q).
\]Therefore
\begin{equation}
Y_{t L_1+(1-t) L_2}(\alpha)\geq tY_{L_1}(\alpha)+(1-t)Y_{L_2}(\alpha)\phantom{.}
\end{equation}
and, in a similar way,
\begin{equation}
X_{t L_1+(1-t) L_2}(\beta)\geq tX_{L_1}(\beta)+(1-t)X_{L_2}(\beta).
\end{equation}
According to equation (\ref{key1}) 
\begin{equation}
f_{t L_1+(1-t) L_2}(x,y)\geq tf_{L_1}(x,y)+(1-t)f_{L_2}(x,y)\ \ \ (0\leq t\leq 1)
\end{equation}
for any $(x,y)\in \mathbb{R}^2$.
\end{Pf}

\begin{Rem} \emph{
Integrating both sides of inequality (3) (or (4)) it follows that  
$$\lambda_2(tL_1+(1-t)L_2)\geq t\lambda_2(L_1)+(1-t)\lambda_2(L_2).$$
If we omit the condition of the common axis parallel bounding box then the statement is false as the following example shows:
$$L_1:=[-3,3]\times [-3,3],\ \ L_2=[-1,1]\times [-1,1],$$
$$\frac{1}{2}L_1+\frac{1}{2}L_2=[-2,2]\times [-2,2]$$
and
$$\lambda_2\left(\frac{1}{2}L_1+\frac{1}{2}L_2\right)=16,\ \ \frac{1}{2}\lambda_2(L_1)+\frac{1}{2}\lambda_2(L_2)=\frac{1}{2}36+\frac{1}{2}4=20.$$
The condition of the common axis parallel bounding box is used as} none of the sets \emph{$\alpha=_1 L_i$ ($i=1,2$)} are empty \emph{because $L_1$ and $L_2\in \mathcal{L}_B$}.
\end{Rem}

Since the generalized conic function is defined by the integral over the focal set $\Phi$ obviously preserves the ordering with respect to the inclusion and the pointwise upper semi-continuity
\[\limsup_{n\to \infty}f_{L_n}(x,y)\leq f_L(x,y)
\]
follows immediately for any $(x,y)\in \mathbb{R}^2$, where the sequence $L_n$ tends to $L\in \mathcal{F}(\mathbb{R}^2)$ with respect to the Hausdorff metric. In the forthcoming section we are going to investigate the continuity properties of the mapping $\Phi$ restricted to the class of compact connected hv-convex planar sets.

\section{The case of compact connected hv-convex sets}

In what follows $\mathcal{F}_{hv}(\mathbb{R}^2)$ denotes the metric space of non-empty compact hv-convex sets in the plane equipped with the Hausdorff metric. Let the rectangle $B$ with parallel sides to the coordinate axes be given as the Cartesian product $[a,b]\times [c,d]$. Recall that $\mathcal{L}_B$ is just the collection of non-empty compact planar sets $L$ for which
\begin{equation}
\label{boxcondition}
\textrm{pr}_1(L)\times \textrm{pr}_2(L)=B.
\end{equation}
The axis parallel bounding box of $L$ is the intersection of all axis parallel boxes that contain $L$. It can be easily seen that any non-empty compact connected hv-convex set $L$ is between the upper and the lower bound functions 
	\begin{equation}
	\label{upperlower}
	h_L(x):=\sup \ \{y\ | \ (x,y)\in L\}\ \ \textrm{and}\ \ g_L(x):=\inf\ \{y\ | \ (x,y)\in L\}
\end{equation}
in the sense that
	\[L=\{(x,y)\ | \ a\leq x\leq b\ \textrm{and}\ g_L(x)\leq y \leq h_L(x)\}.
\]On the other hand
	\[Y_L(x)=h_L(x)-g_L(x)
\]and for any sequence $x_n\to x$
	\[\limsup_{n\to \infty}Y_L(x_n)\leq \limsup_{n\to \infty}h_L(x_n)-\liminf_{n\to \infty}g_L(x_n)\leq h_L(x)-g_L(x)=Y_L(x)
\]which means that the coordinate X-ray function $Y_L$ is upper semi-continuous on the interval $[a,b]$. A similar statement can be formulated in terms of $X_L$. For the general theory of parallel X-rays see [2].

\begin{Lem} The set
\begin{equation}
\label{hvconvexsets}
\mathcal{L}_B^{hv}=\mathcal{L}_B\cap \mathcal{F}_{hv}(\mathbb{R}^2)
\end{equation}
consists of all non-empty compact connected hv-convex sets with axis parallel bounding box $B$. The set
\begin{equation}
\label{hvconvexsets1}
\mathcal{M}_B^{hv}=\mathcal{M}_B\cap \mathcal{F}_{hv}(\mathbb{R}^2)
\end{equation}
consists of all non-empty compact connected hv-convex sets with axis parallel bounding box contained in $B$.
\end{Lem}

\begin{Pf}
For the non-trivial direction let $L$ be a non-empty compact hv-convex set satisfying condition (\ref{boxcondition})
and suppose that $L$ is not connected. This means that there is a non-constant, continuous function $f\colon L\to \{0,1\}$. Since $L$ is hv-convex $f$ must be constant along each horizontal or vertical segments running in $L$. Therefore we can construct a continuous function $\tilde{f}\colon \textrm{pr}_1(L)\to \{0,1\}$ such that $\tilde{f}$ makes the diagram 
$$\textrm{\phantom{valamiv}} L \stackrel{f}{\longmapsto} \{0,1\}$$
$$\textrm{pr}_1 \downarrow \ \ \ \ \ \nearrow $$
$$\textrm{pr}_1(L)$$
commutative. This contradicts to the connectedness of $\textrm{pr}_1(L)=[a,b]$.
\end{Pf}

\begin{Lem} If $L$ is a non-empty compact connected hv-convex set then the outer parallel body $L^{\varepsilon}$ is connected and hv-convex for any $\varepsilon > 0$.
\end{Lem}

\begin{Pf}
To prove that $L^{\varepsilon}$ is hv-convex suppose, in contrary, that it is not true. Without loss of generality we can suppose that the points $Q_1(0,0)$ and $Q_2(0,m)$ belong to the parallel body but the segment joining $Q_1$ and $Q_2$ contains a point $Q(0,y)\notin L^{\varepsilon}$ and $0 < y < m$. Therefore we can choose points $P_1(x_1,y_1)$ and $P_2(x_2,y_2)$ from the intersections of $L$ with the closed disks centered at $Q_1$ and $Q_2$ with radius $\varepsilon$ but $L$ must be disjoint from the closed disk $D$ centered at $Q$ with radius $\varepsilon$. We have that $y_1<y$ because $P_1$ must be under the perpendicular bisector of the segment $Q_1Q$ (otherwise $P_1$ would be closer to $Q$ than to $Q_1$ which is a contradiction). In a similar way $y<y_2$. Therefore $y_1< y < y_2$. On the other hand $x_1\neq x_2$ because the convexity into the vertical direction says that if $x_1=x_2$ then the segment joining $P_1$ and $P_2$ belong to $L$. The parallel body of the segment $P_1P_2$ is a convex set and the points $Q_1$ and $Q_2$ belong to $(P_1P_2)^{\varepsilon}\subset L^{\varepsilon}$. So does the point $Q$ which is a contradiction. Using reflection about the vertical coordinate line if necessary suppose that $x_1<x_2$. We claim that $h_L(x_1)<y$ (for the definition of the upper bound function see (\ref{upperlower})). In opposite case we have
	\[y_1 < y \leq h_L(x_1)
\]
and the vertical line segment joining $P_1(x_1,y_1)$ and $(x_1,h_L(x_1))$ intersects $D$ which is a contradiction. Obviously $h_L(x_2)\geq y_2 > y$. Let us define the number
	\[0\leq s:=\sup \{t\geq 0\ | \ h_L(x)<y\ \ \textrm{for all}\ \ x\in [x_1,x_1+t] \}\leq  x_2-x_1.
\]Then we can choose a sequence $s_n^+\to s$ such that $h_L(s_n^+)\geq y$ and, by the upper semi-continuity of the upper bound function $h_L$, it follows that
	\[y\leq \limsup_{n\to \infty} h_L(s_n^+)\leq h_L(s).
\]We can also choose a sequence $s_n^-\to s$ such that  $g_L(s_n^-)\leq h_L(s_n^-) < y$ and, by the lower semi-continuity of the lower bound function $g_L$, it follows that
	\[y\geq \liminf_{n\to \infty} g_{L}(s_n^-) \geq g_{L}(s).
\]The convexity into the vertical direction gives a contradiction because the (vertical) segment joining $(s,h_L(s))$ with $(s,g_L(s))$ intersects $D$. Since $L$ is connected and hv-convex Lemma 1 implies that
$$\textrm{pr}_1(L)\times \textrm{pr}_2(L)=B,$$
where $B=[a, b]\times [c, d]$ is the axis parallel bounding box of $L$. Therefore
$$\textrm{pr}_1(L^{\varepsilon})\times \textrm{pr}_2(L^{\varepsilon})=[a-\varepsilon, b+\varepsilon]\times [c-\varepsilon, d+\varepsilon]$$
and the connectedness follows by using Lemma 1 again.
\end{Pf}

\begin{Lem}
The limit $L\in \mathcal{F}(\mathbb{R}^2)$ of the sequence $L_n$ of non-empty compact connected hv-convex sets is a connected hv-convex set.
\end{Lem}

\begin{Pf}
We are going to discuss the convexity only into the horizontal direction (the discussion of the vertical direction is similar). Suppose, in contrary, that there exist points $P_1(x_1,y)$, $P_2(x_2,y)\in L$ such that the segment joining $P_1$ with $P_2$ contains a point $Q(x,y)\notin L$. Since the Euclidean distance of $Q$ from $L$ is strictly positive we can choose a positive real number $\varepsilon$ in such a way that $L^{2\varepsilon}$ is disjoint from the closed disk $D$ centered at $Q$ with radius $\varepsilon$. The Hausdorff convergence $L_n\to L$ implies the existence of $n\in \mathbb{N}$ such that $L_n\subset L^{\varepsilon}$ and $L\subset L_n^{\varepsilon}$. Therefore there exist points $R_1,R_2\in L_n$ in the closed disks $D_1$ and $D_1$ centered at $P_1$ and $P_2$ with radius $\varepsilon$, respectively. Since
$L_n^{\varepsilon}\subset (L^{\varepsilon})^{\varepsilon}\subset L^{2\varepsilon}$ we have that $L_n^{\varepsilon}$ is disjoint from $D$.

\begin{figure}
    \centering
        \includegraphics[scale=0.2]{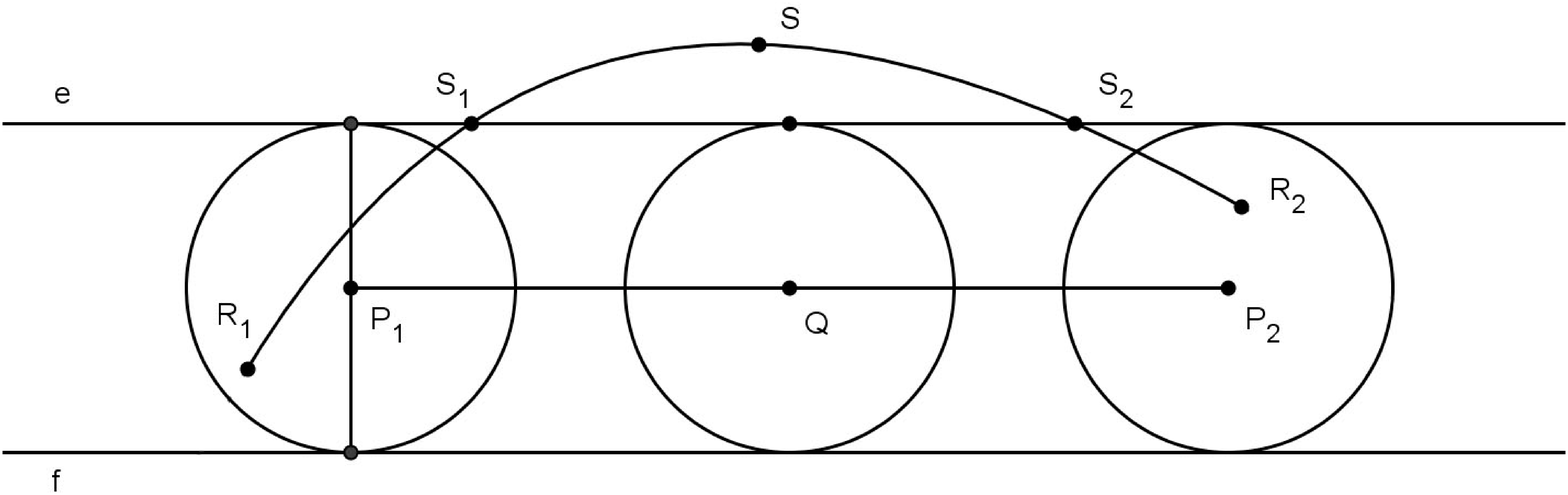}
        \caption{}
\end{figure}

Consider the (common) tangent lines $e$ and $f$ of the disks $D_1$ and $D_2$. They are tangent to $D$ at the same time. Let $H_e$ be the closed half plane bounded by $e$ containing the line $f$ and, in a similar way, $H_f$ denotes the closed half plane bounded by $f$ containing the line $e$. In view of Lemma 2, the parallel body $L_n^{\varepsilon}$ is connected together with its interior containing $L_n$. Therefore $\textrm{int}\ L_n^{\varepsilon}$ is arcwise connected as a connected open subset of the Euclidean plane. There exists a continuous arc $s\subset \textrm{int}\ L_n^{\varepsilon}$ joining $R_1$ and $R_2$ but $L_n^{\varepsilon}$ together with $s$ is disjoint from $D$ (the closed disk centered at $Q$ with radius $\varepsilon$). Taking a point $S\in s$ such that $S\notin H_e\cap H_f$ (suppose, for example, that $S\notin H_e$) we can divide $s$ into the union of continuous arcs $s_1$ and $s_2$ intersecting the line $e$ at the points $S_1$ and $S_2$, respectively. The horizontal segment between $S_1$ and $S_2$ intersects $D$. Since $L_n^{\varepsilon}$ is hv-convex it follows that $D$ and $L_n^{\varepsilon}$ has a common point which is a contradiction.  The connectedness follows easily from Lemma 1 because 
$$B=\textrm{pr}_1(L)\times \textrm{pr}_2(L)$$
is just the limit of the sequence of rectangles $B_n=\textrm{pr}_1(L_n)\times \textrm{pr}_2(L_n)$
with parallel sides to the coordinate axes.
\end{Pf}

\begin{Cor} Both $\mathcal{L}_B^{hv}$ and $\mathcal{M}_B^{hv}$ are compact with respect to the Hausdorff metric. 
\end{Cor}

\begin{Lem}
Let $\varepsilon>0$ be an arbitrary positive real number and consider a finite simple polygonal chain $\mathcal{P}$ in the plane. Then
    \begin{equation}
		\label{poly}
		\lambda_2\left(\mathcal{P}^{\varepsilon}\right)\leq  2l\varepsilon+\varepsilon^2\pi,
		\end{equation}
		where $l$ is the length of $\mathcal{P}$. Especially, if $\mathcal{P}$ is closed then
    \begin{equation}
		\label{polyclosed}
		\lambda_2\left(\mathcal{P}^{\varepsilon}\right)\leq  2l\varepsilon.
\end{equation}
\end{Lem}
\begin{Pf}
The proof is an induction for $m\in\mathbb{N}$, where $m$ denotes the number of segments constituting the chain. If the polygonal chain consists of only one segment then estimation (\ref{poly}) is obviously true (especially we have equality). Suppose that estimation (\ref{poly}) is true for all polygonal chains consisting of $m$ segments and consider a chain $\mathcal{P}_{m+1}$ with vertices $P_1,\ldots,P_{m+1},P_{m+2}$. Taking
	\[\mathcal{P}_{m+1}=\mathcal{P}_m\cup P_{m+1} P_{m+2}
\]we have $\mathcal{P}_{m+1}^{\varepsilon}$ as the union of $\mathcal{P}_m^{\varepsilon}$ and the parallel body of the segment joining $P_{m+1}$ with $P_{m+2}$ (see Figure 2).
\begin{figure}
    \centering
        \includegraphics[scale=0.3]{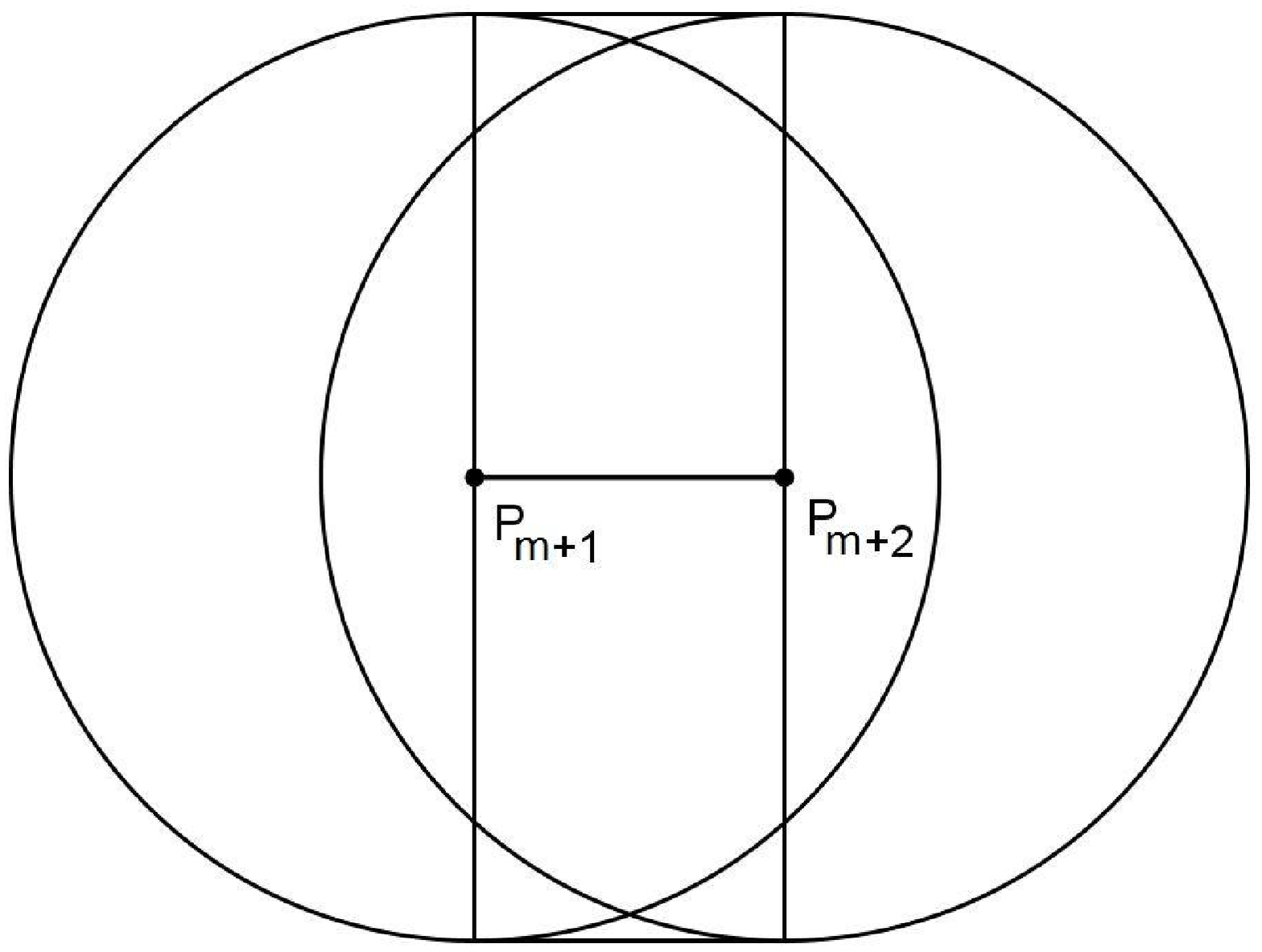}
        \caption{}
\end{figure}
It is clear that the disk around $P_{m+1}$ with radius $\varepsilon$ has a twofold covering. Therefore 
    \[\lambda_2\left(\mathcal{P}_{m+1}^{\varepsilon}\right)\leq \lambda_2\left(\mathcal{P}_{m}^{\varepsilon}\right)+\lambda_2((P_{m+1} P_{m+2})^{\varepsilon})-\varepsilon^2\pi,
\]
where the area of the parallel body of the segment $P_{m+1} P_{m+2}$ is
\[\lambda_2((P_{m+1} P_{m+2})^{\varepsilon})=2\varepsilon\cdot d_2(P_{m+1},P_{m+2})+\varepsilon^2\pi.\]
Therefore
    \[\lambda_2\left(\mathcal{P}_{m+1}^{\varepsilon}\right)\leq \lambda_2\left(\mathcal{P}_{m}^{\varepsilon}\right)+2\varepsilon\cdot d_2(P_{m+1},P_{m+2})
\]
and the inductive hypothesis
gives the estimation
	\[\lambda_2 \left(\mathcal{P}_{m+1}^{\varepsilon}\right)\leq 2l\varepsilon+\varepsilon^2\pi\]
as was to be stated. Estimation (\ref{polyclosed}) follows from (\ref{poly}) by subtracting the areas coming from the twofold coverings of the disks around $P_{m+1}$ and $P_{m+2}=P_1$.
\end{Pf}

\begin{Thm}
Suppose that $L\subset \mathbb{R}^2$ is a non-empty compact connected hv-convex set with axis parallel bounding box $B=[a,b]\times[c,d]$. Then 
\[\lambda_2\left(L^{\varepsilon}\right)-\lambda_2(L)\leq 2k\varepsilon,
\]where $k$ denotes the perimeter of $\ B$.
\end{Thm}

\begin{Pf}
Let $r_n\to 0^+$ be an arbitrary sequence and consider partitions $x_0=a<x_1<\ldots<x_m=b$, $y_0=c<y_1<\ldots<y_m=d$ such that $\diam B^{n}_{ij}< r_n$, where
	\[B^{n}_{ij}:=[x_{i-1},x_{i}]\times [y_{j-1},y_{j}]\quad(i,j=1,2,\ldots,m).
\]
For the minimal covering 
	\[L_n:=\bigcup_{B^{n}_{ij}\cap L\neq \emptyset}B^{n}_{ij}
\] of $L$ as the union of subrectangles having a non-empty intersection with $L$ we have that $H(L_n,L)\leq r_n$, i.e. the sequence $L_n$ tends to $L$ with respect to the Hausdorff metric. On the other hand $\lambda_2(L_n)\to \lambda_2(L)$
because of $L\subset L_n\subset L^{r_n}$. Therefore for any positive real number $\varepsilon > 0$
\[\lambda_2\left(L^{\varepsilon}\right)-\lambda_2(L)=\lambda_2\left(L^{\varepsilon}\right)-\lim_{n\to\infty}\lambda_2(L_n)=\lim_{n\to\infty}\left(\lambda_2\left(L^{\varepsilon}\right)-\lambda_2(L_n)\right)
\]
and thus
	\[\lambda_2\left(L^{\varepsilon}\right)-\lambda_2(L)\leq \lim_{n\to\infty}\left(\lambda_2\left(L_n^{\varepsilon}\right)-\lambda_2(L_n)\right)
\]using that $L\subset L_n$. We claim that $L_n$ is hv-convex. The discussion will be restricted to the convexity into the horizontal direction (the discussion of the vertical direction is similar). Suppose, in contrary, that there exist points $P_1(x_1,y)$, $P_2(x_2,y)\in L_n$ such that the segment joining $P_1$ with $P_2$ contains a point $Q(x,y)\notin L_n$. This means that $P_1$ and $P_2$ are in disjoint subrectangles $B^{\ n}_{i_1 j_1}$ and $B^{\ n}_{i_2 j_2}$, respectively. According to the definition of $L_n$ these subrectangles contain points $R_1$ and $R_2\in L$ but the subrectangle $B^{\ n}_{i_3j_3}$ containing $Q$ must be disjoint from the set $L$. Since $L$ is compact we can choose a positive real number $0< \delta$ such that $L^{\delta}\cap B^{\ n}_{i_3j_3}=\emptyset$. Lemma 2 implies that $L^{\delta}$ is connected together with its interior. Therefore $\textrm{int}\ L^{\delta}$ is arcwise connected as a connected open subset of the Euclidean plane and the points $R_1$, $R_2$ can be joined by a continuous arc in the interior of $L^{\delta}$. The argumentation can be finished in the same way as in the proof of Lemma 3 (see Figure 1 with rectangular domains instead of disks). Since 
$$\textrm{pr}_1(L_n)\times \textrm{pr}_2(L_n)=B$$
it follows by Lemma 2 that it is a connected set: $L_n$ is actually a special kind of parallel body constructed from $L$ by adding rectangular domains instead of disks. Therefore the boundary of $L_n$ is a finite simple closed polygonal chain $\mathcal{P}_n$ in the plane. The lenght of $\mathcal{P}_n$ is the perimeter of the box $B$. Since $L_n^{\varepsilon}\setminus L_n\subset \mathcal{P}_n^{\varepsilon}$
we have by Lemma 4 that
	\[\lambda_2\left(L_n^{\varepsilon}\right)-\lambda_2(L_n)=\lambda_2(L_n^{\varepsilon}\setminus L_n)\leq \lambda_2(\mathcal{P}_n^{\varepsilon})\leq  2k\varepsilon
\]and, consequently,
	\[\lambda_2\left(L^{\varepsilon}\right)-\lambda_2(L)\leq 2k\varepsilon
\]as was to be stated.
\end{Pf}

\begin{Thm}
The mapping $\Phi \colon L\in \mathcal{M}^{hv}_B\to f_L$ is continuous between $\mathcal{M}^{hv}_B$ equip\-ped with the Hausdorff metric and the function space equipped with the norm
	\[\|f_L\|_{\infty, B}:=\sup_{(x,y)\in B} |f_L(x,y)|.
\]
\end{Thm}

\begin{Pf}
Suppose that $L_n\to L$ with respect to the Hausdorff metric, where $L_n$ and $L$ are non-empty compact connected hv-convex sets with axis parallel bounding box contained in $B$. From the definition of the generalized conic function we have that
	\[f_{L_n}(x,y)\leq f_{L^{r_n}}(x,y)=f_L(x,y)+\intop\limits_{L^{r_n}\setminus L} d_1\big{(}(x,y),(\alpha,\beta)\big{)}\, d\alpha d\beta
\]where $r_n:=H(L_n,L)$. The integrand is obviously bounded from above by the following way: since $L^{r_n}\subset B^{r_n}$ and $L\subset B$ we have that
\begin{equation}
d_1\big{(}(x,y),(\alpha,\beta)\big{)}\leq \frac{k}{2}+2r_n
\end{equation}
for any $(x,y)\in B$, where $k$ is the perimeter of $B$: $k/2$ is just the diameter of $B$ with respect to the taxicab norm but the parallel body $L^{r_n}$ allows us two additional steps of lenght $r_n$ into the vertical or the horizontal directions. Therefore
	\[f_{L_n}(x,y)\leq f_L(x,y)+\bigg(\frac{k}{2}+2r_n\bigg) \big(\lambda_2(L^{r_n})-\lambda_2(L)\big) \leq 
\]
	\[f_L(x,y)+\bigg(\frac{k}{2}+2r_n\bigg)2kr_n
\]because of Theorem 2. Conversely
	\[f_{L}(x,y)\leq f_{L_n^{r_n}}(x,y)=f_{L_n}(x,y)+\intop\limits_{L_n^{r_n}\setminus L_n} d_1\big{(}(x,y),(\alpha,\beta)\big{)}\, d\alpha d\beta
\]where $r_n:=H(L_n,L)$.  The integrand is obviously bounded from above by the same way as in (12): 
\begin{equation}
d_1\big{(}(x,y),(\alpha,\beta)\big{)}\leq \frac{k}{2}+2r_n.
\end{equation}
Therefore
	\[f_{L}(x,y)\leq f_{L_n}(x,y)+\bigg(\frac{k}{2}+2r_n\bigg) \big(\lambda_2(L_n^{r_n})-\lambda_2(L_n)\big) \leq 
\]
	\[f_{L_n}(x,y)+\bigg(\frac{k}{2}+2r_n\bigg) 2kr_n
\]because of Theorem 2. These inequalities imply that for any $(x,y)\in B$
\begin{equation}
|f_L(x,y)-f_{L_n}(x,y)|\leq \bigg(\frac{k}{2}+2r_n\bigg)2kr_n,
\end{equation}
where the right hand side is a quadratic polynomial expression of the Hausdorff distance which is independent of the choice of $(x,y)\in B$. In other words the convergence $f_{L_n}\to f_L$ is uniform over the reference set $B$ and the statement of Theorem 3 follows immediately. 
\end{Pf}

\begin{Cor}
The mapping $\Phi \colon L\in \mathcal{M}^{hv}_B\to f_L$ is continuous between $\mathcal{M}^{hv}_B$ equipped with the Hausdorff metric and the function space equip\-ped with the norm
	\[\|f_L\|_{1, B}:=\intop\limits_B |f_L(x,y)|\, dx dy.
\]
\end{Cor}

\begin{Cor}
If $L_n$ is a sequence of non-empty compact connected hv-convex sets tending to the limit $L\in \mathcal{F}(\mathbb{R}^2)$ then for any $(x,y)\in \mathbb{R}^2$
\begin{equation}
\lim_{n\to \infty}f_{L_n}(x,y)=f_L(x,y).
\end{equation}
\end{Cor}

\begin{Pf} Lemma 3 says that $L$ is a non-empty compact connected hv-convex set and we can use Theorem 3 under the choice of a sufficiently large reference set $B$.
\end{Pf}

\begin{Cor}
Let $L_n$ be a sequence of non-empty compact connected hv-convex sets contained in $B$. If $\ f_{L_n}\to f_K$ with respect to the $L_1$-norm or the supremum norm then any convergent subsequence of $L_n$ tends to a set $K'$ having the same coordinate X-rays as $K$ almost everywhere. If $\ K$ is uniquely detemined by the coordinate X-rays then $K'$ is equal to $K$ modulo a set of measure zero.
\end{Cor}

\begin{Pf}
Consider the case of the supremum norm. If $K'$ is the limit of a subsequence $L_{n_k}$ then
	\[\|f_{K'}-f_K\|_{\infty, B}\leq \|f_{K'}-f_{L_{n_k}}\|_{\infty, B}+\|f_{L_{n_k}}-f_K\|_{\infty, B},\]
	where the first term tends to zero in view of Theorem 3. So does the second term because of the condition $f_{L_n}\to f_K$ with respect to the supremum norm.
Taking the limit as $k\to \infty$ 
	\[\|f_{K'}-f_K\|_{\infty, B}=0
\]which means that $f_K(x,y)=f_{K'}(x,y)$ for any $(x,y)\in B$ because of the continuity of the generalized conic functions. Therefore $K$ and $K'$ have the same coordinate X-rays almost everywhere. 
\end{Pf}

\vspace{0.3cm}
As a sequence-free version we can formulate the following theorem of approximation.

\begin{Thm} Suppose that $K\in \mathcal{M}^{hv}_B$. For any $\varepsilon>0$ there exists $\delta_{\infty}>0$ or $\delta_1>0$ such that for any $L\in \mathcal{M}^{hv}_B$
	\[\|f_L-f_K\|_{\infty, B} < \delta_{\infty}\ \ \textrm{or}\ \ \ \|f_L-f_K\|_{1, B} < \delta_1
\]implies that $H(L,K')< \varepsilon$ for some $K'\in \mathcal{M}^{hv}_B$, where $K'$ has the same coordinate X-rays as $K$ almost everywhere. 
\end{Thm}

The theorem says that if $f_L\approx f_K$ with respect to the supremum norm or the $L_1$-norm then $L$ approximates at least one of the sets in $\Phi^{-1}(f_K)$.

\section{The case of compact convex planar bodies: Gardner's problem}

In what follows $\mathcal{F}_{c}(\mathbb{R}^2)$ denotes the metric space of non-empty compact convex sets in the plane equipped with the Hausdorff metric. Let the rectangle $B$ with parallel sides to the coordinate axes be given as the Cartesian product $[a,b]\times [c,d]$. Recall that $\mathcal{L}_B$ is just the collection of non-empty compact planar sets $L$ for which
\begin{equation}
\textrm{pr}_1(L)\times \textrm{pr}_2(L)=B.
\end{equation}
R\r{a}dstr\"om's embedding theorem [4] says that the collection of non-empty compact convex sets (equipped with the Hausdorff metric) can be isometrically embedded into a normed vector space $V$ as a cone. It is a continuous embedding because of the distance preserving property. Therefore 
\begin{equation}
\mathcal{L}_B^c=\mathcal{L}_B\cap \mathcal{F}_{c}(\mathbb{R}^2)
\end{equation}
can be interpreted as a compact convex subset in $V$ (Proposition 1 and Proposition 2) and $\Phi\colon \mathcal{L}_B^c\subset V\to W$ is a bounded order-preserving\footnote{Lattice properties [1] were added to the theory by A. G. Pinsker in 1966.} concave (Theorem 1) and continuous mapping (Theorem 3) into the normed vector space $W$ of continuous functions on $B$ equipped with the supremum or the $L_1$-norm (over B). 

\begin{Def} Let $X$ and $Y$ be Hausdorff topological spaces and consider the mapping $F\colon X \to 2^{Y}$. It is upper semi-continuous at $\ x_0\in X$ if for any open neighbourhood $\mathcal{V}$ of $F(x_0)$ there exists an open neighbourhood $\mathcal{U}$ of $\ x_0$ such that $F(x)\subset \mathcal{V}$ for any $x\in \mathcal{U}$. The mapping $F$ is lower semi-continuous at $\ x_0$ if for any open set $\mathcal{V}$ which intersects $F(x_0)$ there exists an open neighbourhood $\mathcal{U}$ of $x_0$
such that $F(x)\cap \mathcal{V}\neq \emptyset$ for any $x\in \mathcal{U}$. A set-valued mapping which is both upper- and lower semi-continuous is called continuous.
\end{Def}

Since $\Phi\colon L\in \mathcal{L}_B^c \to f_L$ is a continuous mapping defined on a compact metric space its inverse (as a set-valued mapping) is upper semi-continuous at any element $f_K$ of the range: in view of Theorem 4 for any $\varepsilon>0$ there exists $\delta>0$ such that $\|f_L-f_K\|<\delta$ implies that $L\in \big(\Phi^{-1}(f_K)\big)^{\varepsilon}\cap \ \mathcal{L}^c_B$, where
	\[\Phi^{-1}(f_K):=\{K'\in \mathcal{L}^c_B \ | \ f_K(x,y)=f_{K'}(x,y)\ \ \textrm{for all}\ \ (x,y)\in B\}.
\] Under the notation of Definition 3
	\[\mathcal{V}:=\big(\Phi^{-1}(f_K)\big)^{\varepsilon}\cap \ \mathcal{L}^c_B
\]and $\mathcal{U}$ is just the ball with radius $\delta$ around $x_0:=f_K$. The upper semi-continuity establishes an approximating process. In view of Michael's selection theorem the lower semi-continuity is related to the existence of continuous selections. In what follows we prove that the lower semi-continuity of $\Phi^{-1}$ at $f_K$ is equivalent to the determination of $K$ by the coordinate X-rays in the class $\mathcal{F}_{c}^{\circ}(\mathbb{R}^2)$ of non-empty compact convex bodies. A compact convex set is called a body if it has a non-empty interior. The determination of the compact convex body $K$ by the coordinate X-rays means that for any compact convex body $K'$ the relation $f_K=f_{K'}$ implies that $K=K'$. To \emph{characterize those convex bodies that can be determined by two X-rays} is an open problem due to R. J. Gardner [2], Problem 1.1, p. 51. According to the affine nature of the problem we can suppose that the X-ray directions correspond to the coordinate axes without loss of generality. 

\begin{Thm} The body $K\in \mathcal{F}_{c}(\mathbb{R}^2)$  is determined by the coordinate X-rays if and only if the mapping $\ \Phi^{-1}$ is lower semi-continuous  at $f_K$.
\end{Thm}

\begin{Pf} If $K$ is uniquely determined by the coordinate X-rays then Corollary 4 implies immediately the continuity (especially the lower semi-continu\-ity) of the inverse mapping at $f_K$. Conversely, suppose that the inverse mapping is lower semi-continuous. Together with the upper semi-continuity we can conclude that $\Phi^{-1}$  is continuous at $f_K$. Since the set of bodies that can be determined by their coordinate X-rays is dense in $\mathcal{F}_{c}^{\circ}(\mathbb{R}^2)$ [2] [Theorem 1.2.17]  we can choose a sequence $f_{L_n}\to f_K$ such that $\Phi^{-1}(f_{L_n})$ is a singleton. So is $\Phi^{-1}(f_K)$. 
\end{Pf}

\begin{Rem} \emph{The theorem is a way to rephrase the problem of determination in terms of the function $\Phi$. Although it is probably hard to check the continuity of $\Phi^{-1}$ at $f_K$ continuity properties may give answers in terms of algorithms by finding the possible alternatives: let $f_K$ be the input data and consider the optimization problem
$$\textrm{minimize} \ \|f_L-f_K\|\ \ \textrm{subject to}\ \  L\in \mathcal{H},$$ 
where the set $\mathcal{H}$ is the subcollection of non-empty compact connected hv-convex sets which are constituted by the subrectangles belonging to a partition of the reference set \cite{NV4}. It is typically a rectangle $B$ with parallel sides to the coordinate axes such that $K\subset B$.}
\end{Rem}

\section*{Acknowledgements}

The authors would like to thank to the referee for the substantial contribution to the development of the final version of the paper. 

Cs. Vincze was partially supported by the European Union and the European Social Fund through the project Supercomputer, the national virtual lab (grant no.:T\'AMOP-4.2.2.C-11/1/KONV-2012-0010)

Cs. Vincze is supported by the University of Debrecen's internal research project.

\'Abris Nagy has been supported by the Hungarian Academy of Sciences. This research was supported by the European Union and the State of Hungary, co-financed by the European Social Fund in the framework of T\'AMOP-4.2.4.A/ 2-11/1-2012-0001 `National Excellence Program'.

\end{document}